\documentclass[12pt,reqno]{amsart}
\newtheorem{thm}{Theorem}[section]
\newtheorem{defi}{Definition}[section]

\theoremstyle{definition}
\newtheorem{rem}{Remark}[section]
\newtheorem{exm}{Example}[section]
\newcommand{\be}{\begin{equation}}
\newcommand{\ee}{\end{equation}}
\newcommand{\bea}{\begin{eqnarray}}
\newcommand{\eea}{\end{eqnarray}}
\newcommand{\beb}{\begin{eqnarray*}}
\newcommand{\eeb}{\end{eqnarray*}}
\usepackage{amssymb,amsfonts,amsthm,setspace,indentfirst,mathrsfs}
\usepackage{minibox}
\usepackage{enumitem}
\usepackage[utf8x]{inputenc}
\numberwithin{equation}{section}
\begin{document}
%
\title[Curvature properties of a special type of pure radiation metrics]{\bf{Curvature properties of a special type of pure radiation metrics}}
\author[A. A. Shaikh, H. Kundu, M. Ali and Z. Ahsan]{Absos Ali Shaikh$^{*1}$, Haradhan Kundu$^2$, Musavvir Ali$^3$ and Zafar Ahsan$^4$}
\date{\today}
\address{\noindent$^{1,3}$ Department of Mathematics,
\newline Aligarh Muslim University,
\newline Aligarh-202002,
\newline Uttar Pradesh, India}
\email{aask2003@yahoo.co.in, aashaikh@math.buruniv.ac.in}
\email{musavvirali.maths@amu.ac.in}
\address{\noindent$^2$ Department of Mathematics,
\newline University of Burdwan, Golapbag,
\newline Burdwan-713104,
\newline West Bengal, India}
\email{kundu.haradhan@gmail.com}

\address{\noindent$^4$ Faculty of Science and Technology,
\newline University of Islamic Sciences,
\newline Nilai,Malaysia}
\email{zafar.ahsan@rediffmail.com}
\dedicatory{}
\begin{abstract}
A spacetime denotes a pure radiation field if its energy momentum tensor represents a situation in which all the energy is transported in one direction with the speed of light. In 1989, Wils and later in 1997 Ludwig and Edgar studied the physical properties of pure radiation metrics, which are conformally related to a vacuum spacetime. In the present paper we investigate the curvature properties of special type of pure radiation metrics presented by Ludwig and Edgar. It is shown that such a pure radiation spacetime is semisymmetric, Ricci simple, $R$-space by Venzi and its Ricci tensor is Riemann compatible. It is also proved that its conformal curvature 2-forms and Ricci 1-forms are recurrent. We also present a pure radiation type metric and evaluate its curvature properties along with the form of its energy momentum tensor. It is interesting to note that such pure radiation type metric is $Ein(3)$ and 3-quasi-Einstein. We also find out the sufficient conditions for which this metric represents a generalized pp-wave, pure radiation and perfect fluid. Finally we made a comparison between the curvature properties of Ludwig and Edgar's pure radiation metric and pp-wave metrics.
\end{abstract}
%
\maketitle
%

\section{\bf Introduction}\label{intro}
We consider a smooth connected semi-Riemannian manifold $(M^n, g)$ with $n\geq 3$ (this condition is assumed throughout the paper) equipped with the semi-Riemannian metric $g$ with signature $(p, n-p)$ and the Levi-Civita connection $\nabla$. We note that $M$ is Riemannian if $p = 0$ or $n$ and $M$ is Lorentzian if $p = 1$ or $n-1$. Let us consider $R$, $\mathcal R$, $S$, $\mathcal S$ and $\kappa $ respectively be the Riemann-Christoffel curvature tensor of type $(0,4)$, the Riemann-Christoffel curvature tensor of type $(1,3)$, the Ricci tensor of type $(0,2)$, the Ricci tensor of type $(1,1)$ and the scalar curvature of $M$.\\
\indent A spacetime is a connected 4-dimensional Lorentzian manifold and it presents a pure radiation if its energy momentum tensor $T$ is of the form
\be\label{prc}
T = \rho \eta\otimes\eta,
\ee
where $\eta$ is a null vector and $\rho$ is the radiation density. Such a spacetime describes some kinds of field that propagates at the speed of light and it could represent a null electromagnetic field. It could also represent an incoherent beam of photons or some kinds of massless neutrino fields. On the other hand a spacetime presents a perfect fluid spacetime if its energy-momentum tensor is of the form
\be\label{pfc}
T = (\rho + p)\eta\otimes\eta + p g,
\ee
where $\rho$ is the energy density, $p$ is the isotropic pressure and $\eta$ is the four-velocity of the fluid. Thus a pure radiation source can be considered as a limiting case of a pressureless perfect fluid with null four-velocity. For this reason a pure radiation source is sometimes referred as ``null dust''.\\
\indent In \cite{Wils89} Wils investigated homogeneous and conformally Ricci flat solutions of Einstein's field equations for pure radiation case. Later in 1997, Ludwig and Edgar \cite{LE97} obtained exhausted class of conformally Ricci flat pure radiation solutions of Einstein’s field equations. The line element of conformally Ricci flat pure radiation spacetime in $(x^1; x^2; x^3; x^4) = (u; r; x; y)$-coordinates with $x>0$ is given by \cite{LE97}
\beb
ds^2 = \left(-2 x W^{\circ} - 8 \tilde P^2 \frac{r^2}{x^2}\right)du^2 + 2 du dr - \frac{4r}{x} du dx - \frac{1}{8 \tilde P^2}(dx^2 + dy^2),
\eeb
where $\tilde P$ is an arbitrary non-zero constant and $W^{\circ}$ is an arbitrary function of the three non-radial coordinates $u, x$ and $y$.
For simplicity of symbols we write the metric as
\be\label{prm}
ds^2 = \left(x w - p^2 \frac{r^2}{x^2}\right)du^2 + 2 du dr - \frac{4r}{x} du dx - \frac{1}{p^2}(dx^2 + dy^2),
\ee
where $w = w(u,x,y) = -2 W^{\circ}$ and $p = 2\sqrt{2} \tilde P = const.$\\
\indent Now in $(u; r; x; y)$-coordinates with $x > 0$, we consider the following metric
\be\label{gprm}
ds^2 = \left(x w + a \frac{r^2}{x^2}\right) du^2 + 2 du dr + \frac{2 b r}{x} du dx + f (dx^2 + dy^2),
\ee
where $a_1$, $a_2$ are arbitrary non-zero constants and $f$ is a nowhere vanishing function of $x$ and $y$. We note that if $a = -p^2$, $b=-2$ and $f \equiv - \frac{1}{p^2}$, then the metric \eqref{gprm} reduces to pure radiation metric \eqref{prm}. Hence we call the metric \eqref{gprm} as \textit{pure radiation type metric}. Again if $w(u,x,y) = \frac{-2 h(u,x,y)}{x}$, $a = b = 0$ and $f(x,y) = -\frac{1}{2} F(x,y)$, then the metric \eqref{gprm} reduces to generalized pp-wave metric (\cite{RS84}, \cite{SBK17}, \cite{ste03}), given by,
\be\label{gppwm}
ds^2 = -2 h(u,x,y) du^2 + 2 du dr -\frac{1}{2} F (dx^2 + dy^2).
\ee
On the other hand if $w(u,x,y) = \frac{h(u,x,y)}{x}$, $a = b = 0$ and $f \equiv 1$, then the metric \eqref{gprm} reduces to pp-wave metric (\cite{Brin25}, \cite{ste03}), given by,
\be\label{ppwm}
ds^2 = h(u,x,y) du^2 + 2 du dr + (dx^2 + dy^2).
\ee
\indent The physical properties of the pure radiation metric \eqref{prm} are well known and we refer the reader to see \cite{LE97} and \cite{Wils89}. In the literature of differential geometry there are many curvature restricted geometric structures on a semi-Riemannian manifold, such as locally symmetric manifold \cite{Cart26}, semisymmetric manifold (\cite{Cart26}, \cite{Szab82}, \cite{Szab84}, \cite{Szab85}), recurrent manifold (\cite{Ruse46}, \cite{Ruse49a}, \cite{Ruse49b}, \cite{Walk50}), pseudosymmetric manifold (\cite{AD83}, \cite{Desz92} and also references therein) etc. The main object of the present paper is to investigate such kinds of geometric structures admitted by the pure radiation metric \eqref{prm}.
It is noteworthy to mention that the metric \eqref{prm} is neither locally symmetric nor conformally symmetric but semisymmetric and hence Ricci semisymmetric, conformally semisymmetric and projective semisymmetric. It is also shown that the pure radiation metric \eqref{prm} is Ricci simple, weakly Ricci symmetric, weakly cyclic Ricci symmetric, $R$-space by Venzi and its curvature 2-forms, Ricci 1-forms and conformal curvature 2-forms are recurrent. Again the spacetime satisfies the semisymmetric type conditions $C\cdot R =0$, $C\cdot C =0$, $C\cdot S =0$, $Q(S, R) =0$, $Q(S, C) = 0$, $P\cdot R = P\cdot C =0$ and also satisfies the pseudosymmetric type conditions $P\cdot P = -\frac{1}{3}Q(S, P)$.  It is shown that its energy momentum tensor $T$ is semisymmetric and it is Codazzi type (resp., cyclic parallel or covariantly constant) if $w_{33}+w_{44}$ is independent of $x$ and $y$ (resp., constant or zero), where $w_{ij}$ denotes the covariant derivative with respect to $x^i$ and $x^j$.\\
\indent The paper is organized as follows. Section \ref{preli} deals with the preliminaries. In section \ref{com} we compute the components of various tensors of the metric \eqref{prm} and we state the main results on the geometric structures admitted by pure radiation metric \eqref{prm}. Section \ref{gen} deals with the curvature properties of pure radiation type metric \eqref{gprm}. It is shown that such metric is $Ein(3)$ and 3-quasi-Einstein. We also obtain the conditions for which the metric is 2-quasi-Einstein, Ricci generalized pseudosymmetric and manifold of vanishing scalar curvature. Finally, we made a comparison (similarities and dissimilarities) between pure radiation metric and pp-wave metric. It is interesting to mention that both are semisymmetric and weakly Ricci symmetric, but generalized pp-wave metric is Ricci recurrent whereas pure radiation metric is not so.
\section{\bf Curvature Restricted Geometric Structures}\label{preli}
It is wellknown that a curvature restricted geometric structure is a geometric structure on a semi-Riemannian manifold $M$ obtained by imposing a restriction on its curvature tensors by means of covariant derivatives of first order or higher orders. We will now explain some useful notations and definitions of various curvature restricted geometric structures.\\
\indent For two symmetric $(0,2)$-tensors $A$ and $E$, their Kulkarni-Nomizu product $A\wedge E$ is defined as (see e.g. \cite{DGHS11}, \cite{Glog02}):
\begin{eqnarray*}
(A\wedge E)(X_1,X_2,X_3,X_4) &=& A(X_1,X_4)E(X_2,X_3) + A(X_2,X_3)E(X_1,X_4)\\
&-& A(X_1,X_3)E(X_2,X_4) - A(X_2,X_4)E(X_1,X_3),
\end{eqnarray*}
where $X_1,X_2,X_3,X_4 \in \chi(M)$, the Lie algebra of all smooth vector fields on $M$. Throughout the paper we will consider $X, Y, X_1, X_2, \cdots \in \chi(M)$.\\
\indent Again for a symmetric $(0, 2)$-tensor $A$, we get an endomorphism $\mathcal A$ defined by $g(\mathcal AX,Y) = A(X,Y)$. Then its $k$-th level tensor $A^k$ of type $(0,2)$ is given by
$$A^k(X,Y) = A(\mathcal A^{k-1}X,Y),$$
where $\mathcal A^{k-1}$ is the endomorphism corresponding to $A^{k-1}$.\\
\indent In terms of Kulkarni-Nomizu product the conformal curvature tensor $C$, concircular curvature tensor $W$, conharmonic curvature tensor $K$ (\cite{Ishi57}, \cite{YK89}) and the Gaussian curvature tensor $\mathfrak G$ can respectively be expressed as
\beb
C &=& R-\frac{1}{n-2}(g\wedge S) + \frac{r}{2(n-2)(n-1)}(g\wedge g),\\
W &=& R-\frac{r}{2n(n-1)}(g\wedge g),\\
K &=& R-\frac{1}{n-2}(g\wedge S),\\
\mathfrak G &=& \frac{1}{2}(g\wedge g).
\eeb
Again the projective curvature tensor $P$ is given by
$$
P(X_1, X_2, X_3, X_4) = R(X_1, X_2, X_3, X_4) - \frac{1}{n-1}[g(X_1, X_4)S(X_2, X_3)-g(X_2, X_4)S(X_1, X_3)].
$$
\indent For a symmetric $(0,2)$-tensor $A$, $(0,4)$-tensor $D$ and a $(0,k)$-tensor $H$, $k\geq 1$, one can define two $(0,k+2)$-tensors $D\cdot H$ and $Q(A,H)$ respectively as follows (see \cite{DG02}, \cite{DGHS98}, \cite{DH03}, \cite{SDHJK15}, \cite{SK14} and also references therein):
$$D\cdot H(X_1,X_2,\cdots,X_k,X,Y) = -H(\mathcal D(X,Y)X_1,X_2,\cdots,X_k) - \cdots - H(X_1,X_2,\cdots,\mathcal D(X,Y)X_k)$$
and
\beb
Q(A,H)(X_1,X_2, \ldots ,X_k,X,Y) &=& A(X, X_1) H(Y,X_2,\cdots,X_k) + \cdots + A(X, X_k) H(X_1,X_2,\cdots,Y)\\
&-& A(Y, X_1) H(X,X_2,\cdots,X_k) - \cdots - A(Y, X_k) H(X_1,X_2,\cdots,X),
\eeb
where $\mathcal D$ is the corresponding $(1,3)$-tensor of $D$, given by $D(X_1,X_2,X_3,X_4) = g(\mathcal D(X_1,X_2)X_3, X_4).$
\begin{defi}
A semi-Riemannian manifold $M$ is said to be $H$-symmetric (\cite{Cart26}, \cite{Cart46}) if $\nabla H =0$. In particular if $H = R$ (resp., $S$ and $C$), then the manifold is called locally symmetric (resp., Ricci symmetric and conformally symmetric).
\end{defi}
\begin{defi}
A symmetric $(0,2)$-tensor $E$ on $M$ is said to be cyclic parallel (resp, Codazzi type) (see, \cite{DHJKS14}, \cite{Gray78} and references therein) if
\[(\nabla_{X_1} E)(X_2, X_3) = (\nabla_{X_2} E)(X_1, X_3)\]
\[\big(\mbox{resp.,} \ (\nabla_{X_1} E)(X_2, X_3) + (\nabla_{X_2} E)(X_3, X_1) + (\nabla_{X_3} E)(X_1, X_2) = 0\big).\]
\end{defi}
\begin{defi}$($\cite{AD83}, \cite{Cart46}, \cite{Desz92}, \cite{SK14}, \cite{SKppsn}, \cite{SKppsnw}, \cite{Szab82}$)$
A semi-Riemannian manifold $M$ is said to be $H$-semisymmetric type if $D\cdot H = 0$ and it is said to be $H$-pseudosymmetric type if $\left(\sum\limits_{i=1}^k c_i D_i\right)\cdot H = 0$ for some scalars $c_i$'s, where $D$ and each $D_i$, $i=1,\ldots, k$, $(k\ge 2)$, are (0,4) curvature tensors.
\end{defi}
\begin{defi}
A semi-Riemannian manifold $M$ is said to be Einstein if its Ricci tensor is a scalar multiple of the metric tensor $g$. Again $M$ is called quasi-Einstein (resp., 2-quasi-Einstein and 3-quasi-Einstein) if at each point of $M$, rank$(S - \alpha g)\le 1$ (resp., $\le 2$ and $\le 3$) for a scalar $\alpha$. In particular, if $\alpha = 0$, then a quasi-Einstein manifold is called Ricci simple.
\end{defi}
\indent We note that Som-Raychaudhuri metric \cite{SK16srs} and Robinson-Trautman metric \cite{SAArt} are 2-quasi-Einstein whereas G\"{o}del metric \cite{DHJKS14} is Ricci simple.
\begin{defi} (\cite{Bess87}, \cite{SKgrt})
A semi-Riemannian manifold $M$ is said to be $Ein(2)$, $Ein(3)$ and $Ein(4)$ respectively if
$$S^2 + \lambda_1 S + \lambda_2 g = 0,$$
$$S^3 + \lambda_3 S^2 + \lambda_4 S + \lambda_5 g = 0 \ \mbox{and}$$
$$S^4 + \lambda_6 S^3 + \lambda_7 S^2 + \lambda_8 S + \lambda_9 g = 0$$
holds for some scalars $\lambda_i, \, 1\le i \le 9$.
\end{defi}
\begin{defi}\label{def2.6}
Let $D$ be a $(0,4)$-tensor and $E$ be a symmetric $(0, 2)$-tensor on $M$. Then $E$ is said to be $D$-compatible (\cite{DGJPZ13}, \cite{MM12b}, \cite{MM13}) if
\[
D(\mathcal E X_1, X,X_2,X_3) + D(\mathcal E X_2, X,X_3,X_1) + D(\mathcal E X_3, X,X_1,X_2) = 0
\]
holds, where $\mathcal E$ is the endomorphism corresponding to $E$ defined as $g(\mathcal E X_1, X_2) = E(X_1, X_2)$. Again an 1-form $\Pi$ is said to be $D$-compatible if $\Pi\otimes \Pi$ is $D$-compatible.
\end{defi}
\indent Generalizing the concept of recurrent manifold (\cite{Ruse46}, \cite{Ruse49a}, \cite{Ruse49b}, \cite{Walk50}), recently Shaikh et al. \cite{SRK16} introduced the notion of super generalized recurrent manifold along with its characterization and existence by proper example.
\begin{defi}
A semi-Riemannian manifold $M$ is said to be super generalized recurrent manifold (\cite{SK14}, \cite{SKA16}, \cite{SRK16}) if
$$
\nabla R = \Pi \otimes R + \Omega \otimes (S\wedge S) + \Theta \otimes (g\wedge S) + \omega \otimes (g\wedge g)
$$
holds on $\{x\in M: R \neq 0 \mbox{ and any one of } S\wedge S, g\wedge S \mbox{ is non-zero at $x$}\}$ for some 1-forms $\Pi$, $\Omega$, $\Theta$ and $\omega$, called the associated 1-forms. Especially, if $\Omega = \Theta = \omega = 0$ (resp., $\Theta = \omega = 0$ and $\Omega = \omega = 0$), then the manifold is called recurrent (\cite{Ruse46}, \cite{Ruse49a}, \cite{Ruse49b}, \cite{Walk50}) (resp., weakly generalized recurrent (\cite{SAR13}, \cite{SR11}) and hyper generalized recurrent (\cite{SP10}, \cite{SRK15})) manifold.
\end{defi}
\indent Again as a generalization of locally symmetric manifold and recurrent manifold, Tam$\acute{\mbox{a}}$ssy and Binh \cite{TB89} introduced the notion of weakly symmetric manifolds.
\begin{defi}
Let $D$ be a (0, 4)-tensor on a semi-Riemannian manifold $M$. Then $M$ is said to be weakly $D$-symmetric manifold (\cite{TB89}, \cite{SK12}) if
\beb
&&\nabla_X  D(X_1, X_2, X_3, X_4) = \Pi(X) D(X_1, X_2, X_3, X_4) + \Omega(X_1) D(X_1, X_2, X_3, X_4)\\
&& \hspace{2cm} + \overline \Omega(X_2) D(X_1, X, X_3, X_4) + \Theta(X_3) D(X_1, X_2, X, X_4) + \overline \Theta(X_4) D(X_1, X_2, X_3, X)
\eeb
holds $\forall~ X, X_i \in \chi(M)$ $(i =1,2,3,4)$ and some 1-forms $\Pi, \Omega, \overline \Omega, \Theta$ and $\overline \Theta$ on $\{x\in M: R_x \neq 0\}$. Such a manifold is called as weakly $D$-symmetric manifold with solution $(\Pi, \Omega, \overline \Omega, \Theta, \overline \Theta)$. In particular, if the solution is of the form $(2\Pi, \Pi, \Pi, \Pi, \Pi)$, then the manifold is called Chaki $D$-pseudosymmetric manifold \cite{Chak87}. Again if the solution is of the form $(\Pi, 0, 0, 0, 0)$ then the manifold is called $D$-recurrent manifold (\cite{Ruse46}, \cite{Ruse49a}, \cite{Ruse49b}, \cite{Walk50}).
\end{defi}
\begin{defi}
Let $Z$ be a (0, 2)-tensor on a semi-Riemannian manifold $M$. Then $M$ is said to be weakly $Z$-symmetric (\cite{TB93}, \cite{SK12}) if
\beb
(\nabla_X Z)(X_1,X_2)=\Pi(X)\, Z(X_1,X_2) + \Omega(X_1)\, Z(X,X_2) + \Theta(X_2)\, Z(X_1,X)
\eeb
holds $\forall~ X, X_1, X_2 \in \chi(M)$ and some 1-forms $\Pi, \Omega$ and $\Theta$ on $U_{Z} = \{x\in M : Z \neq 0 \ \mbox{ at } x\}$. Such a manifold is called as weakly $Z$-symmetric manifold with solution $(\Pi, \Omega, \Theta)$. Especially, if the solution is of the form $(2\Pi, \Pi, \Pi)$ then the manifold is called Chaki pseudo $Z$-symmetric manifold \cite{Chak88}. Again if the solution is of the form $(\Pi$, $0$, $0)$ then the manifold is called $Z$-recurrent \cite{Patt52}.
\end{defi}
\indent For details about the defining condition of weak symmetry and the interrelation between weak symmetry and Deszcz psudosymmetry, we refer the reader to see \cite{SDHJK15} and also references therein.
\begin{defi}
A Riemannian manifold $M$ is said to be weakly cyclic Ricci symmetric \cite{SJ06} if its Ricci tensor satisfies the condition
\beb
&&(\nabla_X S)(X_1,X_2) + (\nabla_{X_1} S)(X,X_2) + (\nabla_{X_2} S)(X_1,X)\\
&&=\Pi(X)\, S(X_1,X_2) + \Omega(X_1)\, S(X,X_2) + \Theta(X_2)\, S(X_1,X),
\eeb
for three 1-forms $\Pi$, $\Omega$ and $\Theta$ on $M$. Such a manifold is called weakly cyclic Ricci symmetric manifold with solution $(\Pi, \Omega, \Theta)$.
\end{defi}
\indent It is noteworthy to mention that the solution of weakly cyclic Ricci symmetric structure is not always unique.
\begin{defi}
Let $D$ be a $(0,4)$ tensor and $Z$ be a $(0,2)$-tensor on $M$. Then the corresponding curvature 2-forms $\Omega_{(D)l}^m$ (\cite{Bess87}, \cite{LR89}) are said to be recurrent if and only if (\cite{MS12a}, \cite{MS13a}, \cite{MS14})
\beb\label{man}
&&(\nabla_{X_1} D)(X_2,X_3,X,Y)+(\nabla_{X_2} D)(X_3,X_1,X,Y)+(\nabla_{X_3} D)(X_1,X_2,X,Y) =\\
&&\hspace{1in} \Pi(X_1) D(X_2,X_3,X,Y) + \Pi(X_2) D(X_3,X_1,X,Y)+ \Pi(X_3) D(X_1,X_2,X,Y),
\eeb
and the 1-forms $\Lambda_{(Z)l}$ \cite{SKP03} are said to be recurrent if and only if
$$(\nabla_{X_1} Z)(X_2,X) - (\nabla_{X_2} Z)(X_1,X) = \Pi(X_1) Z(X_2,X) - \Pi(X_2) Z(X_1,X)$$
for an 1-form $\Pi$.
\end{defi}
\begin{defi}$($\cite{Prav95}, \cite{SKppsn}, \cite{Venz85}$)$
Let $\mathcal L(M)$ be the vector space formed by all 1-forms $\Theta$ on $M$ satisfying
$$\Theta(X_1)D(X_2,X_3,X_4,X_5)+\Theta(X_2)D(X_3,X_1,X_4,X_5)+\Theta(X_3)D(X_1,X_2,X_4,X_5) = 0,$$
where $D$ is a $(0,4)$-tensor. Then $M$ is said to be a $D$-space by Venzi if $dim [\mathcal L(M)] \ge 1$.
\end{defi}
From definition of recurrency of curvature 2-forms $\Omega_{(R)l}^m$ and second Bianchi identity it is clear that on a semi-Riemannian manifold $\Omega_{(R)l}^m$ are recurrent if and only if it is a $R$ space by Venzi.
\section{\bf Curvature properties of pure radiation metric}\label{com}
The metric tensor of pure radiation metric \eqref{prm} is given by
\beb
g = \left(\begin{array}{cccc}
\left(x w - p^2 \frac{r^2}{x^2}\right) & 1 & \frac{-2r}{x} & 0\\
1 & 0 & 0 & 0\\
\frac{-2r}{x} & 0 & - \frac{1}{p^2} & 0\\
0 & 0 & 0 & - \frac{1}{p^2}\\
\end{array}\right).
\eeb
Then the non-zero components (upto symmetry) of its Riemann-Christoffel curvature tensor $R$, Ricci tensor $S$, scalar curvature $\kappa$, conformal curvature tensor $C$ and projective curvature tensor $P$ are given by
$$R_{1313}=-\frac{w_{33} x}{2}, \ \ R_{1314}=-\frac{w_{34} x}{2}, \ \ R_{1414}=-\frac{w_{44} x}{2},$$
$$S_{11}=-\frac{1}{2} p^2 \left(w_{33}+w_{44}\right) x, \ \ \ \kappa = 0,$$
$$-C_{1313}= C_{1414}=\frac{1}{4} \left(w_{33}-w_{44}\right) x, \ \ C_{1314}=-\frac{w_{34} x}{2},$$
$$P_{1211}=-\frac{1}{6} p^2 \left(w_{33}+w_{44}\right) x, \ \ P_{1311}=\frac{1}{3} p^2 r \left(w_{33}+w_{44}\right),$$
$$P_{1313}=-\frac{1}{6} \left(2 w_{33}-w_{44}\right) x, \ \ -P_{1314}= P_{1341}= -P_{1413}= P_{1431}=\frac{w_{34} x}{2},$$
$$P_{1331}=\frac{w_{33} x}{2}, \ \ P_{1414}=\frac{1}{6} \left(w_{33}-2 w_{44}\right) x, \ \ P_{1441}=\frac{w_{44} x}{2}.$$
Then from above it is easy to check that $R\cdot R = Q(S, R) = R\cdot C = Q(S, C) = 0$.\\
Now the non-zero components (upto symmetry) of $\nabla R$, $\nabla S$, $\nabla C$ are given by
$$R_{1213,1}=\frac{p^2 w_{33}}{2}, \ \ R_{1214,1}=\frac{p^2 w_{34}}{2}, \ \ R_{1313,1}=-\frac{x^2 w_{133}+2 p^2 r w_{33}}{2 x}, \ \ R_{1414,1}=\frac{2 p^2 r w_{44}-x^2 w_{144}}{2 x},$$
$$R_{1313,3}=\frac{1}{2} \left(w_{33}-w_{333} x\right), \ \ R_{1313,4}=-\frac{w_{334} x}{2}, \ \ R_{1314,1}=-\frac{w_{134} x}{2}, \ \ R_{1414,3}=\frac{1}{2} \left(w_{44}-w_{344} x\right),$$
$$R_{1314,3}=\frac{1}{2} \left(w_{34}-w_{334} x\right), \ \ R_{1314,4}=-\frac{w_{344} x}{2}, \ \ R_{1334,1}=\frac{w_{34}}{2}, \ \ R_{1414,4}=-\frac{w_{444} x}{2}, \ \ R_{1434,1}=\frac{w_{44}}{2};$$
$$S_{11,1}=\frac{p^2 \left(-x^2 w_{144}-x^2 w_{133}+2 p^2 r w_{33}+2 p^2 r w_{44}\right)}{2 x},$$
$$S_{11,3}=\frac{1}{2} p^2 \left(-w_{333} x-w_{344} x+w_{33}+w_{44}\right),$$
$$S_{11,4}=-\frac{1}{2} p^2 \left(w_{334}+w_{444}\right) x, \ \ S_{13,1}=\frac{1}{2} p^2 \left(w_{33}+w_{44}\right);$$
$$C_{1213,1}=\frac{1}{4} p^2 \left(w_{33}-w_{44}\right), \ \ C_{1214,1}=\frac{p^2 w_{34}}{2},$$
$$C_{1313,1}=-\frac{-x^2 w_{144}+x^2 w_{133}+2 p^2 r w_{33}-2 p^2 r w_{44}}{4 x},$$
$$C_{1313,3}= -C_{1414,3}=\frac{1}{4} \left(-w_{333} x+w_{344} x+w_{33}-w_{44}\right),$$
$$-C_{1313,4}= C_{1414,4}=\frac{1}{4} \left(w_{334}-w_{444}\right) x, \ \ C_{1314,1}=-\frac{w_{134} x}{2},$$
$$C_{1314,3}=\frac{1}{2} \left(w_{34}-w_{334} x\right), \ \ C_{1314,4}=-\frac{w_{344} x}{2}, \ \ C_{1334,1}=\frac{w_{34}}{2},$$
$$C_{1414,1}=-\frac{x^2 w_{144}-x^2 w_{133}+2 p^2 r w_{33}-2 p^2 r w_{44}}{4 x}, \ \ C_{1434,1}=-\frac{1}{4} \left(w_{33}-w_{44}\right).$$
\indent According to Einstein's field equations, the energy momentum tensor $T$ for zero cosmological constant is related to the Ricci tensor and the metric tensor as
$$T= \frac{c^4}{8\pi G}\left[S-\left(\frac{\kappa}{2}\right)g\right],$$
where $c=$ speed of light in vacuum and $G=$ gravitational constant. Thus the non-zero components of the energy momentum tensor of the pure radiation metric \eqref{prm} is given by:
$$T_{11}=-\frac{c^4 \left(p^2 w_{33} x + p^2 w_{44} x\right)}{16 \pi  G x^2}.$$
Obviously, $T = \rho (\eta\otimes\eta)$, where $\eta = \{1,0,0,0\}$ and the radiation density $\rho = T_{11}$. It is easy to check that $||\eta|| = 0$.\\
Now the non-zero components of covariant derivative of $T$ are given by
$$T_{11,1}=\frac{c^4 p^2 \left(-x^2 w_{144}-x^2 w_{133}+2 p^2 r w_{33}+2 p^2 r w_{44}\right)}{16 \pi  G x}, \ \ T_{13,1}=\frac{c^4 p^2 \left(w_{33}+w_{44}\right)}{16 \pi  G}$$
$$T_{11,3}=\frac{c^4 p^2 \left(-w_{333} x-w_{344} x+w_{33}+w_{44}\right)}{16 \pi  G}, \ \ T_{11,4}=-\frac{c^4 p^2 \left(w_{334}+w_{444}\right) x}{16 \pi  G}.$$
\indent From the value of the local components (presented in Section \ref{com}) of various tensors of the pure radiation metric \eqref{prm}, we can conclude that the pure radiation metric \eqref{prm} fulfills the following curvature restricted geometric structures.
\begin{thm}\label{mainthm}
The pure radiation metric \eqref{prm} possesses the following curvature properties:
\begin{enumerate}[label=(\roman*)]
\item Its Ricci tensor is neither Codazzi type nor cyclic parallel but the scalar curvature is zero and hence $R=W$ and $C =K$.
\item It is a $R$-space (also $C$-space, $P$-space) by Venzi for the associated 1-form $\Pi=\{c,0,0,0\}$, $c$ being arbitrary scalar. Hence the curvature 2-forms $\Omega_{(R)l}^m$ are recurrent for $\Pi$ as the 1-forms of recurrency.
\item It is neither locally symmetric nor conformally symmetric but semisymmetric. Hence it satisfies $R\cdot S=0$, $R\cdot C = 0$ and $R\cdot P = 0$.

\item It satisfies the semisymmetric type condition $C\cdot R=0$ and hence $C\cdot S=0$, $C\cdot C=0$, $C\cdot P=0$ and $P\cdot S = 0$.
\item $R$ or $C$ of the space is not a scalar multiple of $S\wedge S$, but $Q(S,R)=0$, $Q(S,C)=0$. Hence $P\cdot R = 0$ and $P\cdot C = 0$, although $P\cdot \mathcal R\neq 0$ but $P\cdot\mathcal S=0$.
\item It is not Einstein but Ricci simple, since $S = \beta (\eta\otimes\eta)$, where $\beta = - \frac{p^2~x~(w_{33}+w_{44})}{2}$ and $\eta = \{1,0,0,0\}$ (moreover $||\eta|| = 0$ and $\nabla \eta \ne 0$). Hence $S\wedge S=0$ and $S^2=0$.
\item Here $P\cdot P \ne 0$ but $P\cdot P=-\frac{1}{3}Q(S,P)$.
\item If $w_{33}+w_{44}$ is nowhere zero, then the metric is neither recurrent nor Ricci recurrent but Ricci 1-forms are recurrent with associated 1-form
$$\left\{1, 0, \frac{w_{333}+w_{344}}{w_{33}+w_{44}}, \frac{w_{334}+w_{444}}{w_{33}+w_{44}}\right\}.$$
\item If $4w_{34}^2+\left(w_{44}-w_{33}\right)^2$ is nowhere zero, then the metric is not conformally recurrent but conformal 2-forms are recurrent with associated 1-form $\Pi$, given by
$$\Pi_1 = 1, \ \ \ \Pi_2 = 0,$$
$$\Pi_3 = \frac{2 w_{34} \left(w_{334}+w_{444}\right)-\left(w_{44}-w_{33}\right) \left(w_{333}+w_{344}\right)}{4w_{34}^2+\left(w_{44}-w_{33}\right)^2},$$
$$\Pi_4 = \frac{2 w_{34} \left(w_{333}+w_{344}\right)+\left(w_{44}-w_{33}\right) \left(w_{334}+w_{444}\right)}{4w_{34}^2+\left(w_{44}-w_{33}\right)^2}.$$
\item The general form of the compatible tensor for $R$, $C$ and $P$ are respectively given by
$$\left(
\begin{array}{cccc}
 a_{(1,1)} & a_{(1,2)} & a_{(1,3)} & a_{(1,4)} \\
 a_{(2,1)} & 0 & 0 & 0 \\
 a_{(3,1)} & 0 & a_{(3,3)} & a_{(3,4)} \\
 a_{(4,1)} & 0 & a_{(4,3)} & a_{(3,3)}+\frac{w_{44} a_{(3,4)}}{w_{34}}-\frac{w_{33} a_{(4,3)}}{w_{34}}
\end{array}
\right),$$
$$\left(
\begin{array}{cccc}
 a_{(1,1)} & a_{(1,2)} & a_{(1,3)} & a_{(1,4)} \\
 a_{(2,1)} & 0 & 0 & 0 \\
 a_{(3,1)} & 0 & a_{(3,3)} & a_{(3,4)} \\
 a_{(4,1)} & 0 & -\frac{2 w_{34} a_{(3,3)}}{w_{44}-w_{33}}-a_{(3,4)}+\frac{2 w_{34}
   a_{(4,4)}}{w_{44}-w_{33}} & a_{(4,4)}
\end{array}
\right) \ \mbox{and}$$
$$\left(
\begin{array}{cccc}
 a_{(1,1)} & a_{(1,2)} & a_{(1,3)} & a_{(1,4)} \\
 a_{(2,1)} & 0 & 0 & 0 \\
 a_{(3,1)} & 0 & -\frac{\left(2 w_{44}-w_{33}\right) a_{(3,4)}}{3 w_{34}}-\frac{\left(w_{44}-2 w_{33}\right)
   a_{(4,3)}}{3 w_{34}}+a_{(4,4)} & a_{(3,4)} \\
 a_{(4,1)} & 0 & a_{(4,3)} & a_{(4,4)}
\end{array}
\right),$$
where $a_{(i,j)}$ are arbitrary scalars.
\item Ricci tensor of this spacetime is not Codazzi type but is compatible for $R$, $C$ and $P$.
\item If $w_{33}+w_{44}$ is nowhere vanishing, then the metric is weakly Ricci symmetric with infinitely many solutions $(\Pi, \Omega, \Theta)$, given by
$$\Pi = \left\{\Pi_1, 0, \frac{w_{333} x+w_{344} x-w_{33}-w_{44}}{\left(w_{33}+w_{44}\right) x}, \frac{w_{334}+w_{444}}{w_{33}+w_{44}}\right\},$$
$$\Omega = \left\{\Omega_1, 0, -\frac{1}{x}, 0\right\} \ \ \mbox{and}$$
$$\Theta = \left\{\frac{x^2 w_{144}x^2 w_{133}-2 p^2 r w_{33}-2 p^2 rw_{44}}{\left(w_{33}+w_{44}\right) x^2}-\Pi_1-\Omega_1, 0, -\frac{1}{x}, 0\right\},$$
where $\Pi_1$ and $\Omega_1$ are arbitrary scalars.
\item If $w_{33}+w_{44}$ is nowhere vanishing, then the metric is weakly cyclic Ricci symmetric with infinitely many solutions $(\Pi, \Omega, \Theta)$, given by
$$\Pi = \left\{\Pi_1, 0, \frac{\left(w_{333}+w_{344}\right) x-3 w_{33}-3 w_{44}}{\left(w_{33}+w_{44}\right) x}, \frac{w_{334}+w_{444}}{w_{33}+w_{44}}\right\},$$
$$\Omega = \left\{\Omega_1, 0, \frac{\left(w_{333}+w_{344}\right) x-3 w_{33}-3 w_{44}}{\left(w_{33}+w_{44}\right) x}, \frac{w_{334}+w_{444}}{w_{33}+w_{44}}\right\} \ \ \mbox{and}$$
\beb
&&\Theta = \left\{\frac{3 \left(x^2 \left(w_{144}+w_{133}\right)-2 p^2 r \left(w_{33}+w_{44}\right)\right)}{\left(w_{33}+w_{44}\right) x^2}-\Pi_1-\Omega_1, 0,\right.\\
&&\hspace{2.4in} \left.\frac{\left(w_{333}+w_{344}\right) x-3 w_{33}-3 w_{44}}{\left(w_{33}+w_{44}\right) x}, \frac{w_{334}+w_{444}}{w_{33}+w_{44}}\right\},
\eeb
where $\Pi_1$ and $\Omega_1$ are arbitrary scalars.
\item It is not weakly symmetric for $R$, $C$, $P$, $W$ and $K$ and hence not Chaki pseudosymmetric for $R$, $C$ or $P$.
\item $div R \ne 0$, $div C \ne 0$, $div P \ne 0$.
\end{enumerate}
\end{thm}
%
\indent Now from the values of the non-zero components of $\nabla T$, we get
\be\label{cst}
\left.\begin{array}{l}
T_{11,1}+T_{11,1} + T_{11,1}=\frac{3 c^4 p^2 \left(-x^2 w_{144}-x^2 w_{133}+2 p^2 r w_{33}+2 p^2 r w_{44}\right)}{16 \pi  G x},\\
T_{11,3} + T_{13,1} + T_{31,1}=\frac{c^4 p^2 \left(-w_{333} x-w_{344} x+3 w_{33}+3 w_{44}\right)}{16 \pi  G},\\
T_{11,4} + T_{14,1} + T_{41,1}=-\frac{c^4 p^2 \left(w_{334}+w_{444}\right) x}{16 \pi  G}
\end{array}\right\}\ee
and
\be\label{codt}
\left.\begin{array}{l}
T_{13,1} - T_{11,3} = \frac{c^4 p^2 \left(w_{333}+w_{344}\right) x}{16 \pi  G},\\
T_{14,1} - T_{11,4} = \frac{c^4 p^2 \left(w_{334}+w_{444}\right) x}{16 \pi  G}.
\end{array}\right\}\ee
\indent Again since $R\cdot S = 0$ and $T$ is a linear combination of $S$ and $g$, so $R\cdot T = 0$. Hence we can state the following:
\begin{thm}
The energy-momentum tensor $T$ of the pure radiation metric \eqref{prm} is\\
(i) semisymmertric i.e., $R\cdot T = 0$,\\
(ii) Codazzi type if $w_{33}+w_{44}$ is independent of $x$ and $y$,\\
(iii) cyclic parallel if $w_{33}+w_{44}$ is independent of $u$, $x$ and $y$,\\
(iv) covariantly constant if $w_{33}+w_{44} = 0$.
\end{thm}
\section{\bf Curvature properties of pure radiation type metric}\label{gen}
We now consider the pure radiation type metric \eqref{gprm}. Its metric tensor is given by
\beb
g = \left(\begin{array}{cccc}
\left(x w + a \frac{r^2}{x^2}\right) & 1 & \frac{b r}{x} & 0\\
1 & 0 & 0 & 0\\
\frac{b r}{x} & 0 & f & 0\\
0 & 0 & 0 & f\\
\end{array}\right),
\eeb
where $a, b$ are arbitrary non-zero constants and $w = w(u,x,y)$, $f = f(x,y)$ are nowhere vanishing functions.
Then the non-zero components (upto symmetry) of $R$, $S$ and $\kappa$ are given by
$$R_{1212}=-\frac{4 a f-b^2}{4 f x^2}, \ \ R_{1213}=\frac{r \left(8 a f+b^3\right)}{4 f x^3}, \ \ R_{1334}=-\frac{b^2 f_4 r}{4 f x^2},$$
\beb
R_{1313}=-\frac{1}{4 f x^4}&&\left[a b^2 f r^2+2 a b f_3 r^2 x+2 a f_3 r^2 x+12 a f r^2-b^4 r^2+b^2 f w x^3\right.\\
&&\left.+2 b f w_3 x^4+2 b f w x^3-f_3 w_3 x^5+f_4 w_4 x^5+2 f w_{33} x^5-f_3 w x^4+4 f w_3 x^4\right],
\eeb
$$R_{1314}=-\frac{2 a b f_4 r^2+2 a f_4 r^2+b f w_4 x^3-f_4 w_3 x^4-f_3 w_4 x^4+2 f w_{34} x^4-f_4 w x^3+2 f w_4 x^3}{4 f x^3},$$
$$R_{1323}=-\frac{b \left(b f+f_3 x+2 f\right)}{4 f x^2}, \ \ R_{1324}= R_{1423}=-\frac{b f_4}{4 f x},$$
$$R_{1414}=-\frac{-2 a b f_3 r^2-2 a f_3 r^2+f_3 w_3 x^4-f_4 w_4 x^4+2 f w_{44} x^4+f_3 w x^3}{4 f x^3},$$
$$R_{1424}=\frac{b f_3}{4 f x}, \ \ R_{1434}=\frac{b^2 f_3 r}{4 f x^2}, \ \ R_{3434}=-\frac{-f_3^2-f_4^2+f f_{33}+f f_{44}}{2 f},$$
$$S_{11}=\frac{2 a^2 f r^2-3 a b^2 r^2-8 a b r^2+2 a f w x^3-6 a r^2-b(b w + w_3 x + w) x^3-(w_{33} + w_{44}) x^5-2 w_3 x^4}{-2 f x^4},$$
$$S_{12}=-\frac{2 a f-b^2-b}{2 f x^2}, \ \ S_{44}=-\frac{b f f_3+f_3^2 x+f_4^2 x-f f_{33} x-f f_{44} x}{2 f^2 x}, \ \ S_{34}=\frac{b f_4}{2 f x},$$
$$S_{13}=\frac{r \left(4 a f+b^3+b^2\right)}{2 f x^3}, \ \ S_{33}=\frac{b^2 f^2+2 b f^2+b f_3 f x+f_{33} f x^2+f_{44} f x^2-f_3^2 x^2-f_4^2 x^2}{2 f^2 x^2},$$
$$\kappa = \frac{-4 a f^3+b (3 b+4) f^2+2 \left(f_{33}+f_{44}\right) f x^2-2 \left(f_3^2+f_4^2\right) x^2}{2 f^3 x^2}.$$
\indent Again for zero cosmological constant, the non-zero components (upto symmetry) of the energy momentum tensor are given by
\beb
T_{11}=\frac{c^4}{32 \pi  f^3 G x^4} &&\left[3 a b^2 f^2 r^2+12 a b f^2 r^2+12 a f^2 r^2+2 a f_3^2 r^2 x^2+2 a f_4^2 r^2 x^2-2 a f f_{33} r^2 x^2\right.\\
&&-2 a f f_{44} r^2 x^2-b^2 f^2 w x^3+2 b f^2 w_3 x^4-2 b f^2 w x^3+2 f^2 w_{33} x^5+2 f^2 w_{44} x^5\\
&&\left. +4 f^2 w_3 x^4+2 f_3^2 w x^5+2 f_4^2 w x^5-2 f f_{33} w x^5-2 f f_{44} w x^5\right],
\eeb
$$T_{12}=-\frac{c^4 \left(b^2 f^2+2 b f^2+2 f_{33} f x^2+2 f_{44} f x^2-2 f_3^2 x^2-2 f_4^2 x^2\right)}{32 \pi  f^3 G x^2},$$
$$T_{13}=\frac{c^4 r \left(4 a b f^3+8 a f^3+b^3 \left(-f^2\right)-2 b^2 f^2+2 b f_3^2 x^2+2 b f_4^2 x^2-2 b f f_{33} x^2-2 b f f_{44} x^2\right)}{32 \pi  f^3 G x^3},$$
$$T_{33}=\frac{c^4 \left(4 a f^2+b^2 (-f)+2 b f_3 x\right)}{32 \pi  f G x^2}, \ \ T_{34}=\frac{b c^4 f_4}{16 \pi  f G x}, \ \ T_{44}=\frac{c^4 \left(4 a f^2-3 b^2 f-2 b f_3 x-4 b f\right)}{32 \pi  f G x^2}.$$
From above we have the following:
\begin{thm}
The pure radiation type metric \eqref{gprm} has the following curvature properties:\\
(i) It is a 3-quasi-Einstein manifold, since $S- \frac{b+b^2-2 a f}{2 f x^2} g$ is of rank 3.\\
(ii) For zero cosmological constant, its energy momentum tensor $T$ is of the form
$$T = \alpha g + \beta_{_1} (e_1\otimes e_1) + \beta_{_3} (e_3\otimes e_3) + \beta_{_4} (e_4\otimes e_4) +  \sigma_{_1} (e_1\otimes e_3 + e_3\otimes e_1) + \sigma_{_2} (e_3\otimes e_4 + e_4\otimes e_3),$$
where
$$e_1 = \{1,0,0,0\}, \ \ e_3 = \{0,0,1,0\}, \ \ e_4 = \{0,0,0,1\},$$
$$\alpha =-\frac{c^4 \left[b (b+2) f^2+2 \left(f_{33}+f_{44}\right) f x^2-2 \left(f_3^2+f_4^2\right) x^2\right]}{32 \pi  f^3 G x^2},$$
$$\beta_{_1} = T_{11} - \alpha g_{11}, \ \ \beta_{_3} = T_{33} - \alpha g_{33}, \ \ \beta_{_4} = T_{44} - \alpha g_{44}, \ \ 
\sigma_{_1} = T_{13} - \alpha g_{13}, \ \ \sigma_{_2} = T_{34} - \alpha g_{34}.$$
(iii) It is an $Ein(4)$ manifold, such that
\beb
S^4 &=&\frac{\left(-2 a f+b^2+b\right)^2}{16 f^8 x^7}\left[b (b+2) f^3 \left(\left(f_{33}+f_{44}\right) x-b f_3\right)\right.\\
&&+f^2 x  \left(\left(f_{33}+f_{44}\right){}^2 x^2-2 b (b+1) \left(f_3^2+f_4^2\right)\right)\\
&&\left.-2 \left(f_3^2+f_4^2\right) \left(f_{33}+f_{44}\right) f x^3+\left(f_3^2+f_4^2\right)^2 x^3\right] g\\
&-& \frac{\left(-2 a f+b^2+b\right)}{8 f^7 x^6} \left[f^4 \left(b^2 (b+1) (b+2)-4 a \left(f_{33}+f_{44}\right) x^2\right)\right.\\
&&+2 f^3 x \left(x \left(2 a f_4^2+b (2 b+3) \left(f_{33}+f_{44}\right)\right)+2 a f_3^2 x-(b+2) b^2 f_3\right)\\
&&-2 a b (b+2) f^5+2 f^2 x^2 \left(\left(f_{33}+f_{44}\right)^2 x^2-3 b (b+1) \left(f_3^2+f_4^2\right)\right)\\
&&\left.-4 \left(f_3^2+f_4^2\right)\left(f_{33}+f_{44}\right) f x^4+2 \left(f_3^2+f_4^2\right)^2 x^4\right] S\\
&+& \frac{1}{4 f^6 x^4}\left[4 a^2 f^6+f^4 \left(b^2 (b+1) (3 b+5)-8 a \left(f_{33}+f_{44}\right) x^2\right)\right.\\
&&+f^3 x \left(x \left(8 a f_4^2+b (5 b+6) \left(f_{33}+f_{44}\right)\right)+8 a f_3^2 x-(b+2) b^2 f_3\right)\\
&&-4 a b (2 b+3) f^5+f^2 x^2 \left(\left(f_{33}+f_{44}\right){}^2 x^2-6 b (b+1) \left(f_3^2+f_4^2\right)\right)\\
&&\left.-2 \left(f_3^2+f_4^2\right) \left(f_{33}+f_{44}\right) f x^4+\left(f_3^2+f_4^2\right)^2 x^4\right] S^2\\
&+& \frac{1}{2 f^3 x^2}\left[f^2 (4 a f-b (3 b+4))+2 \left(f_3^2+f_4^2-f \left(f_{33}+f_{44}\right)\right) x^2\right] S^3,
\eeb
\end{thm}
\begin{rem}
Since for zero cosmological constant $S = \frac{\kappa}{2}g +\frac{8\pi G}{c^4} T$, hence from above theorem we can state that the Ricci tensor of \eqref{gprm} is of the form
\beb
S &=& (S_{12}) g + (S_{11} - S_{12} g_{11}) (e_1\otimes e_1) + (S_{33} - S_{12} g_{33}) (e_3\otimes e_3) + (S_{44} - S_{12} g_{44}) (e_4\otimes e_4)\\
&& +  (S_{13} - S_{12} g_{13}) (e_1\otimes e_3 + e_3\otimes e_1) + (S_{34} - S_{12} g_{34}) (e_3\otimes e_4 + e_4\otimes e_3).\eeb
\end{rem}
\indent We can easilly check that $||e_1||=0$ and the nonzero components of $\nabla e_1$ are given by
$$\nabla_{_1} e_1 = \frac{a r}{x^2}, \ \ \nabla_{_3} e_1 = \frac{b}{2 x}.$$
We know that a spacetime is called generalized pp-wave (\cite{RS84}, \cite{SBK17}) if there exists a covariantly constant null vector field. Hence we can state the following:
\begin{thm}
The pure radiation type metric \eqref{gprm} represents generalized pp-wave if $a = b = 0$.
\end{thm}
%
\indent Now a spacetime is 2-quasi-Einstein if Rank$(S-\alpha g) = 0$. Hence the metric \eqref{gprm} becomes 2-quasi-Einstein if one of the following condition holds
$$(i) \ (S_{34} - S_{12} g_{34}) = (S_{44} - S_{12} g_{44}) = 0,$$
$$(ii) \ (S_{13} - S_{12} g_{13}) = (S_{33} - S_{12} g_{33}) = (S_{34} - S_{12} g_{34}) = 0,$$
$$(iii) \ (S_{11} - S_{12} g_{11}) = (S_{13} - S_{12} g_{13}) = 0.$$
Simplifying the above conditions we can state the following:
\begin{thm}
The pure radiation type metric \eqref{gprm} becomes 2-quasi-Einstein if any one the following condition holds
$$\begin{array}{l}
(i) \ b = 2 a f^3-\left(f_3^2+f_4^2-f \left(f_{33}+f_{44}\right)\right) x^2=0,\\
(ii) \ f_4 = 2 a f^3+x f \left(x f_{33}-b f_3\right)-b (b+1) f^2-x^2 f_3^2 = 0,\\
(iii) \ a = b = -f_3^2+f_4^2+f \left(f_{33}+f_{44}\right) = 0,\\
(iv) \ f_4 = a = b f_3+\frac{b f}{x}+x f_{33}-\frac{x f_3^2}{f} = 0,\\
(v) \ f_4 = b-2 = \frac{2 a f^2}{x}+x f_{33}-2 f_3-\frac{x f_3^2}{f}-\frac{2 f}{x}=0,\\
(vi) \ a = (b+2) w_3+w_{33} x + w_{44} x = 0,\\
(vii) \ b-2 = w_{33}+w_{44} = 0.
\end{array}$$
\end{thm}
\begin{exm}
If we consider the metric \eqref{gprm} with $f(x, y) = \frac{e^{\frac{1}{3} x^3}}{x^{2/3}}$, $b = -2$ and $a = 0$, then \eqref{gprm} becomes a 2-quasi-Einstein manifold.
\end{exm}
\indent Now a spacetime is perfect fluid if it satisfies \eqref{pfc}. Hence the metric \eqref{gprm} represents perfect fluid if one of the following condition holds
$$\mbox{(i) } \beta_{_1} = \beta_{_3} = \sigma_{_1} = \sigma_{_2} = 0,$$
$$\mbox{(ii) } \beta_{_1} = \beta_{_4} = \sigma_{_1} = \sigma_{_2} = 0,$$
$$\mbox{(iii) } \beta_{_3} = \beta_{_4} = \sigma_{_1} = \sigma_{_2} = 0.$$
Simplifying the above conditions we can state the following:
\begin{thm}
The pure radiation type metric \eqref{gprm} represents perfect fluid if any one the following condition holds
$$\begin{array}{l}
(i) \ a = b = 2 w_3 + x w_{33} + x w_{44} = f_3^2+f_4^2-f \left(f_{33}+f_{44}\right) = 0,\\
(ii) \ b+2 = f_4 = w_{33}+w_{44} = 2 a f^3 - x^2 f_3^2 + x f \left(x f_{33} - 2 f_3\right) - 2 f^2=0,\\
(iii) \ a = f_4 = (b+2) w_3+w_{33} x+w_{44} x = x f \left(b f_3+x f_{33}\right)+b f^2-x^2 f_3^2 = 0,\\
(iv) \ a = f_4 = x f \left(b f'+x f''\right)+b f^2-x^2 f'^2 = x f \left(b f'-x f''\right)+b (b+1) f^2+x^2 f'^2 = 0,\\
(v) \ b+2 = f-\frac{1}{a} = 0.
\end{array}$$
\end{thm}
\begin{exm}
Consider the metric \eqref{gprm}, where $w(u, x, y) = u x y$, $f(x, y) = \frac{e^{\frac{1}{3} x^3}}{x^{2/3}}$, $b = -2$ and $a = 0$, then \eqref{gprm} represents a perfect fluid spacetime. The energy momentum tensor of this metric can be expressed as
$$T = -\frac{c^4 e^{-\frac{x^3}{3}} \left(3 x^3+1\right)}{24 \pi  G x^{4/3}} g + \frac{c^4 \left(3 x^3-2\right)}{12 \pi  G x^2} (e_4 \otimes e_4).$$
Moreover in this case the metric is quasi-Einstein and $Ein(2)$.
\end{exm}
\indent Again a spacetime is a pure radiation spacetime if it satisfies \eqref{prc}. Now $||e_1||=0$ and hence the metric \eqref{gprm} represents perfect fluid if $\alpha = \beta_{_3} = \beta_{_4} = \sigma_{_1} = \sigma_{_2} = 0$. Simplifying these conditions we can state the following:
\begin{thm}
The pure radiation type metric \eqref{gprm} represents pure radiation if $f \equiv \frac{1}{a}$ and $b = -2$.
\end{thm}
Again from above we can easily calculate (large but straightforward) the components of $R\cdot R$, $Q(g,R)$ and $Q(S,R)$. Then we have the following:
\begin{thm}
The pure radiation type metric \eqref{gprm} is\\
(i) Ricci generalized pseudosymmetric ($R\cdot R = Q(S,R)$) if $b = -2$  and $f$ is constant. And in this case the metric is $R$-space by Venzi for the associated 1-form $\{1,0,0,0\}$,\\
(ii) a manifold of vanishing scalar curvature if $f\equiv \frac{b(3 b + 4)}{4 a}$,\\
(iii) semisymmetric if $b = -2$ and $f \equiv \frac{1}{a}$.
\end{thm}
\indent From \eqref{prm}, \eqref{gprm} and \eqref{gppwm}, we see that both the pure radiation metric and the pp-wave metric are special cases of the metric \eqref{gprm}. We refer the reader to see \cite{MS16}, \cite{SK-pp} and also references therein for recent works on pp-wave metric. We now draw a comparison (similarities and dissimilarities) between the curvature properties of pure radiation metric and pp-wave metric.\\
\textbf{A. Similarities:}
\begin{enumerate}
	\item Both the metrics are of vanishing scalar curvature,
	\item both are $R$-space by Venzi as well as $C$-space by Venzi,
	\item both are semisymmetric and semisymmetric due to conformal curvature tensor,
	\item for both the metrics $Q(S,R) = Q(S,C) = 0$,
	\item both the metrics are Ricci simple,
	\item Ricci tensors of both metrics are Riemann compatible as well as conformal compatible,
	\item Ricci 1-forms $\Lambda_{(Z)l}$ of both metrics are recurrent,
	\item conformal 2-forms $\Omega_{(C)l}^m$ of both metrics are recurrent,
	\item for both the metrics $P\cdot R = 0$ but $P\cdot\mathcal R \ne 0$.
	\item both are weakly Ricci symmetric and hence weakly cyclic Ricci symmetric,
	\item both the metrics satisfy $P\cdot P = -\frac{1}{3}Q(S, P)$,
	\item the energy-momentum tensors of both the metrics are semisymmetric.
\end{enumerate}
\textbf{B. Dissimilarities:}
\begin{enumerate}
	\item For zero cosmological constant, the energy-momentum tensors of both the metrics are of rank one, but the associated 1-form for radiation metric is null, and for pp-wave metric it is null as well as covariantly constant,
	\item pp-wave metric is Ricci recurrent but pure radiation metric is not so,
	\item for pp-wave metric energy-momentum tensor is cyclic parallel if and only if it is parallel but this fact is not true for pure radiation metric.
\end{enumerate}


\end{document}